\documentclass[graybox,natbib]{svmult}
\usepackage{mathptmx}
\usepackage{helvet}
\usepackage{courier}
\usepackage{type1cm}
\usepackage{makeidx}
\usepackage{graphicx}
\usepackage{url}
\usepackage{multicol}
\usepackage{chapterbib}
\usepackage[bottom]{footmisc}
\usepackage{enumerate}

\usepackage[english]{babel}
\usepackage{amsmath,amssymb}

\def\IR{\mathbb{R}}

\def\E{{\rm E}}

\DeclareMathOperator{\conv}{\,conv}

\def\ddd#1#2#3{\mathop{{\it D}^{#3}}\nolimits(#1|#2)}

\def\dddd#1#2#3{\mathop{{\it D}^{#3}}\nolimits
({ #1}|{ #2}^1,\dots,{ #2}^n)}
\def\co{{\it co}}
\def\conv{{\it co}}

\def\med{\mbox{\rm med}}
\def\sumi{\sum_{i=1}^n}
\def\V{{\rm vol}}

\makeindex

\begin{document}

\title*{Depth statistics}
\author{Karl Mosler}
% Use \authorrunning{Short Title} for an abbreviated version of author names if the original ones are too long
\institute{Karl Mosler \at Universit\"at zu K\"oln, Albertus-Magnus-Platz, 50923 K\"oln, Germany, \\ \email{mosler@statistik.uni-koeln.de}}
\maketitle
\section{Introduction}\label{sec1}
In 1975 John Tukey proposed a multivariate median which is the `deepest' point in a given data cloud in $\IR^d$ \citep{Tukey75}.
 In measuring the depth of an arbitrary point $z$ with respect to the data, \cite{DonohoG92}
considered hyperplanes through $z$ and determined its `depth' by the smallest portion of data that are separated by such a hyperplane.
Since then, this idea has proved extremely fruitful.
A rich statistical methodology has developed that is based on data depth and, more general, nonparametric depth statistics.
General notions of data depth have been introduced as well as many special ones. These notions vary regarding their computability and robustness and their sensitivity to reflect asymmetric shapes of the data. According to their different properties they fit to particular applications.
The upper level sets of a depth statistic provide a family of set-valued statistics, named \emph{depth-trimmed} or \emph{central regions}. They describe the distribution regarding its location, scale and shape. The most central region serves as a \emph{median}.
%%; see also Chapter \ref{Oja}.
The notion of depth has been extended from data clouds, that is empirical distributions, to general probability distributions on $\mathbb{R}^d$, thus allowing for laws of large numbers and consistency results. It has also been extended from $d$-variate data to data in functional spaces.
The present chapter surveys the theory and methodology of depth statistics.

Recent reviews on data depth are given in \cite{Cascos09} and \cite{Serfling06}. \cite{LiuSS06} collects theoretical as well as applied work.  More on the theory of depth functions and many details are found in \cite{ZuoS00a} and the monograph by \cite{Mosler02a}.

The depth of a data point is reversely related to its \emph{outlyingness}, and the depth-trimmed regions can be seen as multivariate set-valued \emph{quantiles}. To illustrate the notions we consider bivariate data from the EU-27 countries regarding unemployment rate and general government debt in percent of the GDP (Table 1). In what follows we are interested which countries belong to a central, rather homogeneous group and which have to be regarded as, in some sense, outlying.

\begin{table}
\begin{center}
  \begin{tabular}{|l|r|r||l|r|r|}
  \hline
  Country & Debt \% & Unempl. \% & Country & Debt \% & Unempl. \% \\
  \hline
  Belgium & 98.0 & 7.2 & Luxembourg & 18.2 & 4.9 \\
  Bulgaria & 16.3 & 11.3 & Hungary & 80.6 & 10.9 \\
  Czech Republic & 41.2 & 6.7 &  Malta & 72.0 & 6.5 \\
  Denmark & 46.5 & 7.6 & Netherlands & 65.2 & 4.4 \\
  Germany & 81.2 & 5.9 & Austria & 72.2 & 4.2 \\
  Estonia & 6.0 & 12.5 & Poland & 56.3 & 9.7 \\
  Ireland & 108.2 & 14.4 & Portugal & 107.8 & 12.9 \\
  Greece & 165.3 & 17.7 & Romania & 33.3 & 7.4 \\
  Spain & 68.5 & 21.7 & Slovenia & 47.6 & 8.2 \\
  France & 85.8 & 9.6 & Slovakia & 43.3 & 13.6 \\
  Italy & 120.1 & 8.4 & Finland & 48.6 & 7.8 \\
  Cyprus & 71.6 & 7.9 & Sweden & 38.4 & 7.5 \\
  Latvia & 42.6 & 16.2 & United Kingdom & 85.7 & 8.0 \\
  Lithuania & 38.5 & 15.4 & & & \\
  \hline
  \end{tabular}
\end{center}
\caption{
General government gross debt (\% of GDP) and
unemployment rate of the EU-27 countries in 2011 (Source: EUROSTAT)}.
\end{table}

Overview: Section \ref{sec2} introduces general depth statistics and the notions related to it. In Section \ref{sec3} various depths for $d$-variate data are surveyed: multivariate depths based on distances, weighted means, halfspaces or simplices.
Section \ref{sec4} provides an approach to depth for functional data, while Section \ref{sec5} treats computational issues.
Section \ref{sec6} concludes with remarks on applications.

\section{Basic concepts}\label{sec2}
In this section the basic concepts of depth statistics are introduced, together with several related notions.
First we provide a general notion of depth functions, which relies on a set of desirable properties; then a few variants of the properties are discussed (Section 2.1). A depth function induces an outlyingness function and a family of central regions (Section 2.2). Further, a stochastic ordering and a probability metric are generated (Section 2.3).

\subsection{Postulates on a depth statistic}\label{subsec2.1}
Let $E$ be a Banach space, ${\cal B}$ its Borel sets in $E$, and ${\cal P}$ a set of  probability distributions on ${\cal B}$.
To start with and in the spirit of Tukey's approach to data analysis, we may regard ${\cal P}$ as the class of empirical distributions giving equal probabilities $\frac 1n$ to $n$, not necessarily different, data points in $E=\mathbb{R}^d$.

A \emph{depth function} is a function $D:E\times {\cal P} \to [0,1], \; (z,P) \mapsto D(z|P),$ that satisfies the restrictions (or `postulates') \textbf{D1} to \textbf{D5} given below.
For easier notation we write $D(z|X)$ in place of $D(z|P)$, where $X$ is an arbitrary random variable distributed as $P$.
For $z\in E$, $P\in {\cal P}$, and any random variable $X$ having distribution $P$ it holds:
\begin{itemize}
  \item \textbf{D1 Translation invariant:} $D(z+b|X+b)=D(z|X)$ for all $b\in E$\,.
  \item \textbf{D2 Linear invariant:} $D(Az|AX)=D(z|X)$ for every bijective linear transformation $A: E\to E$\,.
  \item \textbf{D3 Null at infinity:} $\lim_{\left\|z\right\|\rightarrow\infty}D(z|X)=0$\,.
  \item \textbf{D4 Monotone on rays:} If a point $z^*$ has maximal depth, that is $D(z^*|X)=\max_{z\in E}D(z|X)$\,,
  then for any $r$ in the unit sphere of $E$ the function $\alpha\mapsto D(z^*+\alpha r|X)$ decreases, in the weak sense, with $\alpha>0$\,.
  \item \textbf{D5 Upper semicontinuous:} The upper level sets  $D_\alpha(X) =  \{z\in E : D(z|X) \geq \alpha\}$ are closed for all $\alpha$\,.
\end{itemize}
\textbf{D1} and \textbf{D2} state that a depth function is \emph{affine invariant}. \textbf{D3} and \textbf{D4} mean that the level sets $D_\alpha$, $\alpha>0$, are bounded and starshaped about $z^*$. If there is a point of maximum depth, this depth will w.l.o.g.\ be set to 1. \textbf{D5} is a useful technical restriction.
An immediate consequence of restriction \textbf{D4} is the following:
\begin{proposition}\label{prop:symmetric}
   If $X$ is centrally symmetric distributed about some $z^*\in E$, then any depth function $D(\cdot|X)$ is maximal at $z^*$.
\end{proposition}
Recall that $X$ is \emph{centrally symmetric} distributed about $z^*$ if the distributions of $X-z^*$ and $z^*-X$ coincide.

Our definition of a depth function differs slightly from that given in \cite{Liu90} and  \cite{ZuoS00a}.
The main difference between these postulates and ours is that they additionally postulate Proposition \ref{prop:symmetric} to be true and that they do not require upper semicontinuity \textbf{D5}.

\textbf{D4} states that the upper level set $D_\alpha(x^1,\dots, x^n)$ are starshaped with respect to $z^*$.
If a depth function, in place of \textbf{D4}, meets the restriction
\begin{itemize}
  \item \textbf{D4con:} $D(\cdot|X)$ is a \textbf{quasiconcave} function, that is, its upper level sets $D_\alpha(X)$  are convex for all $\alpha >0$\,,
\end{itemize}
the depth is mentioned as a \emph{convex depth}.
Obviously, as a convex set is starshaped with respect to each of its points,
$\textbf{D4con}$ implies $\textbf{D4}$.
In certain settings the restriction \textbf{D2} is weakened to
\begin{itemize}
  \item \textbf{D2iso:} $D(Az|AX)=D(z|X)$ for every \textbf{isometric linear} transformation $A: E\to E$\,.
\end{itemize}
Then, in case $E=\IR^d$, $D$ is called an \emph{orthogonal invariant depth} in contrast to an \emph{affine invariant} depth when \textbf{D2} holds. Alternatively, sometimes \textbf{D2} is attenuated to \textbf{scale invariance},
\begin{itemize}
  \item \textbf{D2sca:} $D(\lambda z|\lambda X)=D(z|X)$ for all $\lambda > 0$\,.
\end{itemize}

\subsection{Central regions and outliers}\label{subsec2.2}
For given $P$ and $0\le \alpha\le 1$ the level sets $D_\alpha(P)$ form a nested family of \textit{depth-trimmed} or \textit{central regions}.
  The innermost region arises at some $\alpha_{max}\le 1$, which in general depends on $P$. $D_{\alpha_{max}}(P)$ is the set of \textit{deepest points}.
\textbf{D1} and \textbf{D2} say that the family of central regions is affine equivariant.
  %%\textbf{D3} means that for any $\alpha>0$ the region $D_\alpha(X)$ is bounded. \textbf{D4} states the starshapedness of each $D_\alpha(X)$ with respect to %%$z^*$.
Central regions describe a distribution $X$  with respect to location, dispersion, and shape. This has many applications in multivariate data analysis.
On the other hand, given a nested family $\{C_\alpha(P)\}_{\alpha\in [0,1]}$ of set-valued statistics, defined on ${\cal P}$, that are convex, bounded and closed, the function $D$,
\begin{equation}\label{eqregiondepth}
    D(z|P) = \sup\{\alpha : z\in C_\alpha(P)\}\,,\quad z\in E, \;  P \in {\cal P},
\end{equation}
satisfies \textbf{D1} to \textbf{D5} and \textbf{D4con}, hence is a convex depth function.

A depth function $D$ orders data by their degree of centrality. Given a sample, it provides a center-outward \emph{order statistic}. The depth induces an \emph{outlyingness function} $\IR^d\to [0,\infty[$ by
\[ Out(z|X) =\frac 1 {D(z|X)} -1\,,\]
which is zero at the center and infinite at infinity. In turn, $D(z|X) = (1+Out(z|X))^{-1}$.
Points outside a central region $D_\alpha$ have outlyingness greater than $1/\alpha -1$; they can be regarded as \emph{outliers} of a specified level $\alpha$.

\subsection{Depth lifts, stochastic orderings, and metrics}\label{subsec2.3}
Assume $\alpha_{max}=1$ for $P \in {\cal P}$. By adding a real dimension to the central regions $D_\alpha(P), \alpha\in [0,1]$, we construct a set, which will be mentioned as the \emph{depth lift},
\begin{equation}\label{eqdepthlift}
\widehat D (P)= \{(\alpha, y) \in [0,1]\times E : y=\alpha x,\, x\in D_\alpha (P), \, \alpha\in [0,1]\}\,.
\end{equation}

The depth lift gives rise to an \emph{ordering} of probability distributions in ${\cal P}$: \, $P \prec_D Q$ if
\begin{equation}\label{eqdepthordering}
 \widehat D(P) \subset \widehat D(Q)\,.
\end{equation}
The restriction $\widehat D(P) \subset \widehat D(Q)$ is equivalent to $D_\alpha(P) \subset D_\alpha(Q)$ for all $\alpha$.
Thus, $P \prec_D Q$ means that each central set of $Q$ is larger than the respective central set of $P$. In this sense, $Q$ is \emph{more dispersed} than $P$.
The depth ordering is antisymmetric, hence an \emph{order}, if and only if the family of central regions completely characterizes the underlying probability. Otherwise it is a preorder only.
Finally, the depth $D$ introduces a \emph{probability semi-metric} on ${\cal P}$ by the Hausdorff distance of depth lifts,
\begin{equation}\label{eqdepthnorm}
 \delta_D(P,Q) = \delta_H(\widehat D(P), \widehat D(Q))\,.
\end{equation}
Recall that the \emph{Hausdorff distance} $\delta_H(C_1,C_2)$ of two
compact sets $C_1$ and $C_2$ is the smallest $\varepsilon$ such that $C_1$ plus
the $\varepsilon$-ball includes $C_2$ and vice versa. Again, the semi-metric is a metric iff the central regions characterize the probability.

\section{Multivariate depth functions}\label{sec3}

Originally and in most existing applications depth statistics are used with data in Euclidean space. Multivariate depth statistics are particularly suited to analyze non-gaussian or, more general, non-elliptical distributions in $\mathbb{R}^d$. Without loss of generality, we consider distributions of full dimension $d$, that is, whose convex hull of support, $\conv(P)$, has affine dimension $d$.

A random vector $X$ in $\IR^d$ has a \emph{spherical distribution} if $AX$ is distributed as $X$ for every orthogonal matrix $A$. It has an \emph{elliptical distribution} if $X=a+BY$ for some $a\in \IR^d$, $B\in \IR^{d\times d}$, and spherically distributed $Y$; then we write $X\sim {\rm Ell}(a,BB',\varphi)$, where $\varphi$ is the radial distribution of $Y$. Actually, on an elliptical distribution $P={\rm Ell}(a,BB',\varphi)$, any depth function $D(\cdot,P)$ satisfying ${\bf D1}$ and ${\bf D2}$ has parallel elliptical level sets $D_\alpha(P)$, that is, level sets of a quadratic form with
\emph{scatter matrix} $BB'$. Consequently, all affine invariant depth functions are essentially equivalent if the distribution is elliptical.
Moreover, if $P$ is elliptical and has a unimodal Lebesgue-density $f_P$, the density level sets have the same elliptical shape, and the density is a transformation of the depth, i.e., a function $\varphi$ exists such that $f_P(z)= \varphi(D(z|P)$ for all $z\in \IR^d$.
Similarly, on a spherical distribution, any depth satisfying postulates ${\bf D1}$ and ${\bf D2iso}$ has analogous properties.

In the following, we consider three principal approaches to define a multivariate depth statistic. The first approach is based on distances from properly defined central points or on volumes (Section 3.1), the second on certain L-statistics (\emph{viz.} decreasingly weighted means of order statistics;  Section 3.2),
the third on simplices and halfspaces in $\IR^d$  (Section 3.3). The three approaches have different consequences on the depths' ability to reflect asymmetries of the distribution, on their robustness to possible outliers, and on their computability with higher-dimensional data.

Figures \ref{fig1} to \ref{fig4} below exhibit bivariate central regions for several depths and equidistant $\alpha$. The data consist of the unemployment rate (in \%) and the GDP share of public debt for the countries of the European Union in 2011.

 Most of the multivariate depths considered are convex and affine invariant, some exhibit spherical invariance only.
 Some are continuous in the point $z$ or in the distribution $P$ (regarding weak convergence), others are not.
 They differ in the shape of the depth lift and whether it uniquely determines the underlying distribution.
%%, hence the depth ordering is antisymmetric.
A basic dispersion ordering of multivariate probability distributions serving as a benchmark is the \emph{dilation order}, which says that $Y$ spreads out more than $X$ if $\E[\varphi(X)] \le \E[\varphi(Y)]$ holds for every convex $\varphi: \IR^d\to\IR$; see, e.g.\ \cite{Mosler02a}. It is interesting whether or not a particular depth ordering is concordant with the dilation order.

\subsection{Depths based on distances}\label{subsec3.1}

The outlyingness of a point, and hence its depth, can be measured by a distance from a properly chosen center of the distribution.
In the following notions this is done with different distances and centers.

{\bf $L_2$-depth.}
The $L_2$-depth, $D^{L_2}$, is based on the mean
outlyingness of a point, as measured by the $L_2$ distance,
\begin{equation}\label{L2depth}
\ddd{z}{X}{L_2} = \left( 1 + \E||z-X||\right)^{-1} \, .
\end{equation}
It holds $\alpha_{max}=1$. The depth lift is $\widehat D^{L_2}(X)=\{(\alpha,z) : \E||z- \alpha X||\le 1 -\alpha\}$ and convex.
For an empirical distribution on points $x^i, i=1,\dots,n,$ we obtain
\begin{equation}\label{L2depthsample}
\dddd{z}{x}{L_2} = \left( 1 + \frac 1n
              \sumi ||z-x^i||\right)^{-1}\,.
\end{equation}
Obviously, the $L_2$-depth vanishes at infinity (\textbf{D3}), and is maximum at the {\em spatial
median} of $X$, i.e., at the point $z\in \IR^d$ that minimizes $\E||z-X||$.
If the distribution is centrally symmetric, the center
is the spatial median, hence the maximum is attained at
the center. Monotonicity with respect to the deepest point (\textbf{D4})
as well as convexity and compactness of the central regions (\textbf{D4con}, \textbf{D5}) derive
immediately from the triangle inequality. Further, the $L_2$-depth depends continuously on $z$.
The $L_2$-depth converges also in the probability distribution: For a uniformly integrable and weakly convergent sequence $P_n\to P$ it holds
$\lim_n D(z|P_n)= D(z|P)$.

However, the ordering induced by the $L_2$-depth is no sensible ordering of dispersion, since
the $L_2$-depth contradicts the dilation order.
As $||z-x||$ is convex in $x$, the expectation
$\E||z-X||$ increases with a dilation of $P$.
Hence (\ref{L2depth}) decreases (!) with a dilation.

The $L_2$-depth is invariant against rigid Euclidean
motions (\textbf{D1}, \textbf{D2iso}),  but not affine invariant.
An affine invariant version is constructed as follows:
Given a positive definite $d\times d$ matrix $M$,
consider the {\em M-norm},
\begin{equation}\label{Mahnorm}
||z||_M= \sqrt{z' M^{-1} z}, \quad z\in \IR^d\, .
\end{equation}
Let $S_{X}$ be a positive definite
$d\times d$ matrix that depends continuously (in weak convergence) on
the distribution and measures the dispersion of
$X$ in an affine equivariant way. The latter means that
\begin{equation}\label{affequdisp}
S_{X A + b}= {A}S_{X}{A'} \quad
\mbox{holds for any matrix}\; A \; \mbox{of full rank and any} \; b .
\end{equation}
Then an {\em affine invariant $L_2$-depth}\/ is given by
\begin{equation}\label{affinL2depth}
         \left(1 + \E||z-X||_{S_{X}}\right)^{-1} \, .
\end{equation}
Besides invariance, it has the same properties as the $L_2$-depth.
A simple choice for $S_{X}$
is the covariance matrix ${\Sigma}_{X}$ of $X$ \citep{ZuoS00a}.
Note that the covariance matrix is positive definite, as
the convex hull of the support, $\conv(P)$, is assumed to have full dimension.
More robust choices for $S_{X}$ are the {\it minimum volume ellipsoid}\/ (MVE) or  the
{\it minimum covariance determinant}\/ (MCD) estimators; see \citet{RousseeuwL87},
\citet{LopuhaaR91}.
%%, and Chapter \ref{Rousseeuw} of this volume.

\begin{figure}
\includegraphics[scale=0.39]{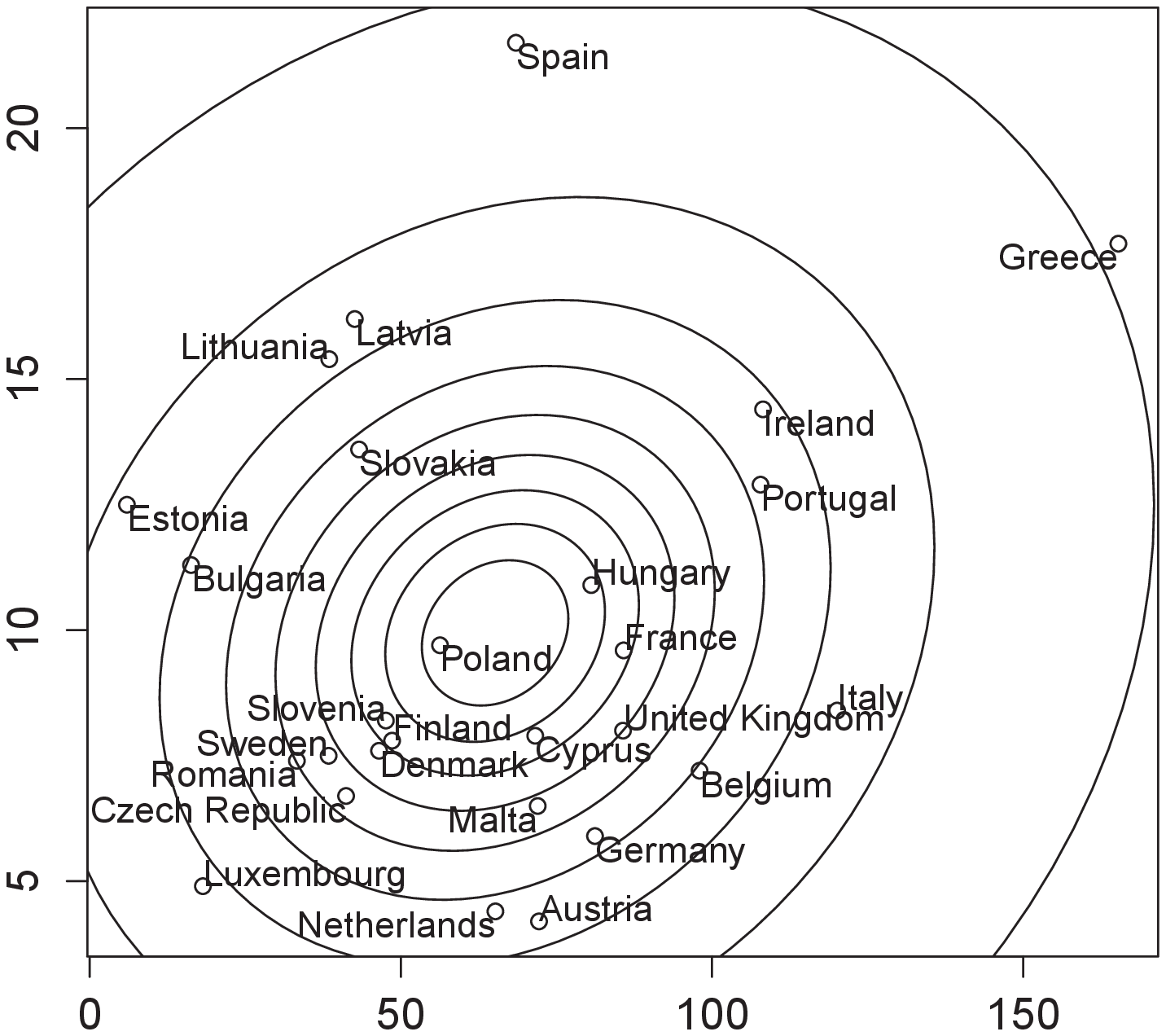}
\quad \includegraphics[scale=0.39]{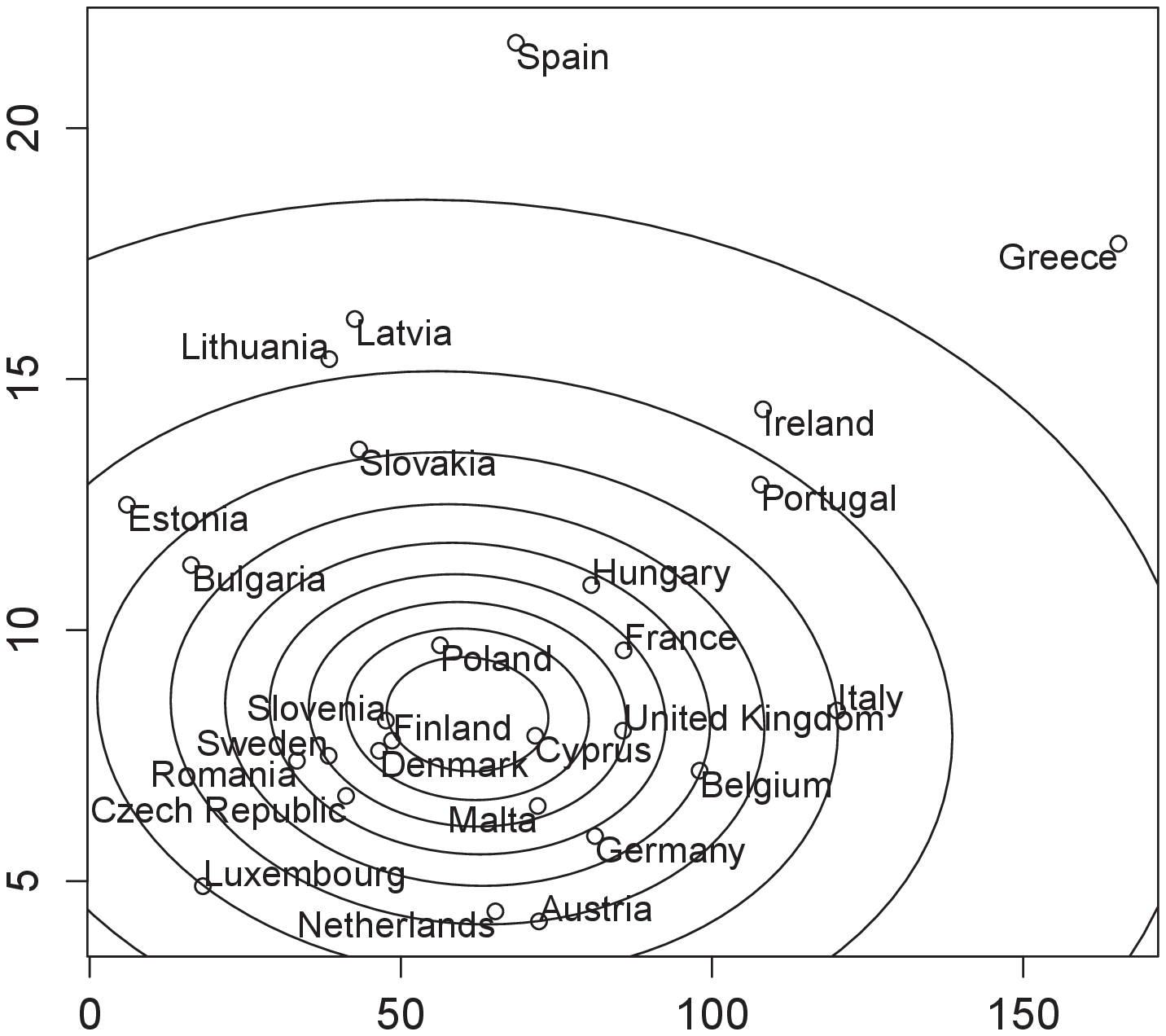}
\caption{Governmental debt  ($x$-axis) and unemployment rate ($y$-axis); Mahalanobis regions (moment, left; MCD, right) with $\alpha = 0.1 (0.1),\dots, 0.9$.}
\label{fig1}
\end{figure}

{\bf Mahalanobis depths.}
Let $c_{X}$ be a vector that measures the location of $X$
in a continuous and affine equivariant way and, as before,
$S_{X}$ be a matrix that satisfies
(\ref{affequdisp}) and depends continuously on the distribution.
Based on the estimates $c_{X}$ and $S_{X}$
a simple depth statistic is constructed, the {\em generalized Mahalanobis depth}, given by
\begin{equation}\label{gMahdepth}
\ddd{z}{X}{Mah} =
         \left( 1 + ||z-c_{X}||^2_{S_{X}}\right)^{-1} \,.
\end{equation}
Obviously, (\ref{gMahdepth}) satisfies \textbf{D1} to \textbf{D5} and \textbf{D4con}, taking its unique maximum at $c_X$.
The depth lift is the convex set
$\widehat D^{Mah}(X)=\{(\alpha,z) : ||z - \alpha c_X||^2_{S_X} \le \alpha^2(\alpha - 1)\}$, and the
central regions are ellipsoids around  $c_X$. The generalized Mahalanobis depth is continuous on $z$ and $P$.
In particular, with $c_{X}=\E[X]$ and $S_{X}=\Sigma_X$ the {\em (moment) Mahalanobis depth}
is obtained,
\begin{equation}\label{pMahdepth}
\ddd{z}{X}{mMah} =
         \Bigl( 1 + (z-\E[X])' {\Sigma}_{X}^{-1}
              (z-\E[X]) \Bigr)^{-1} \, .
\end{equation}
Its sample version is
\begin{equation}\label{sMahdepth}
\dddd{z}{x}{mMah} =
         \Bigl( 1 + (z-\overline x)' {\widehat\Sigma}_{x}^{-1}
              (z-\overline x) \Bigr)^{-1} \, ,
\end{equation}
where $\overline x$ is the mean vector and  ${\widehat\Sigma}_{X}$ is the empirical covariance matrix.
It is easily seen that the $\alpha$-central set
of a sample from $P$ converges almost surely to the
$\alpha$-central set of $P$, for any $\alpha$. Figure \ref{fig1} shows Mahalanobis regions for the debt-unemployment data,
employing two choices of the matrix $S_X$, namely the usual moment estimate $\Sigma_X$ and the robust
MCD estimate. As it is seen from the Figure, these region depend heavily on the choice of $S_X$. Hungary, e.g., is rather central (having depth greater than 0.8) with the moment Mahalanobis depth, while it is much more outlying (having depth below 0.5) with the MCD version.

Concerning uniqueness, the Mahalanobis depth
fails in identifying the underlying distribution.
As only the first two moments are used, any two distributions which
have the same first two moments cannot be distinguished by
their Mahalanobis depth functions. Similarly,
the {generalized Mahalanobis depth} does not determine the distribution.
However, within the family of nondegenerate
$d$-variate normal distributions
or, more general, within any affine family of nondegenerate
$d$-variate distributions having finite second moments,
a single contour set of the Mahalanobis depth suffices to identify
the distribution.

\begin{figure}
\includegraphics[scale=0.39]{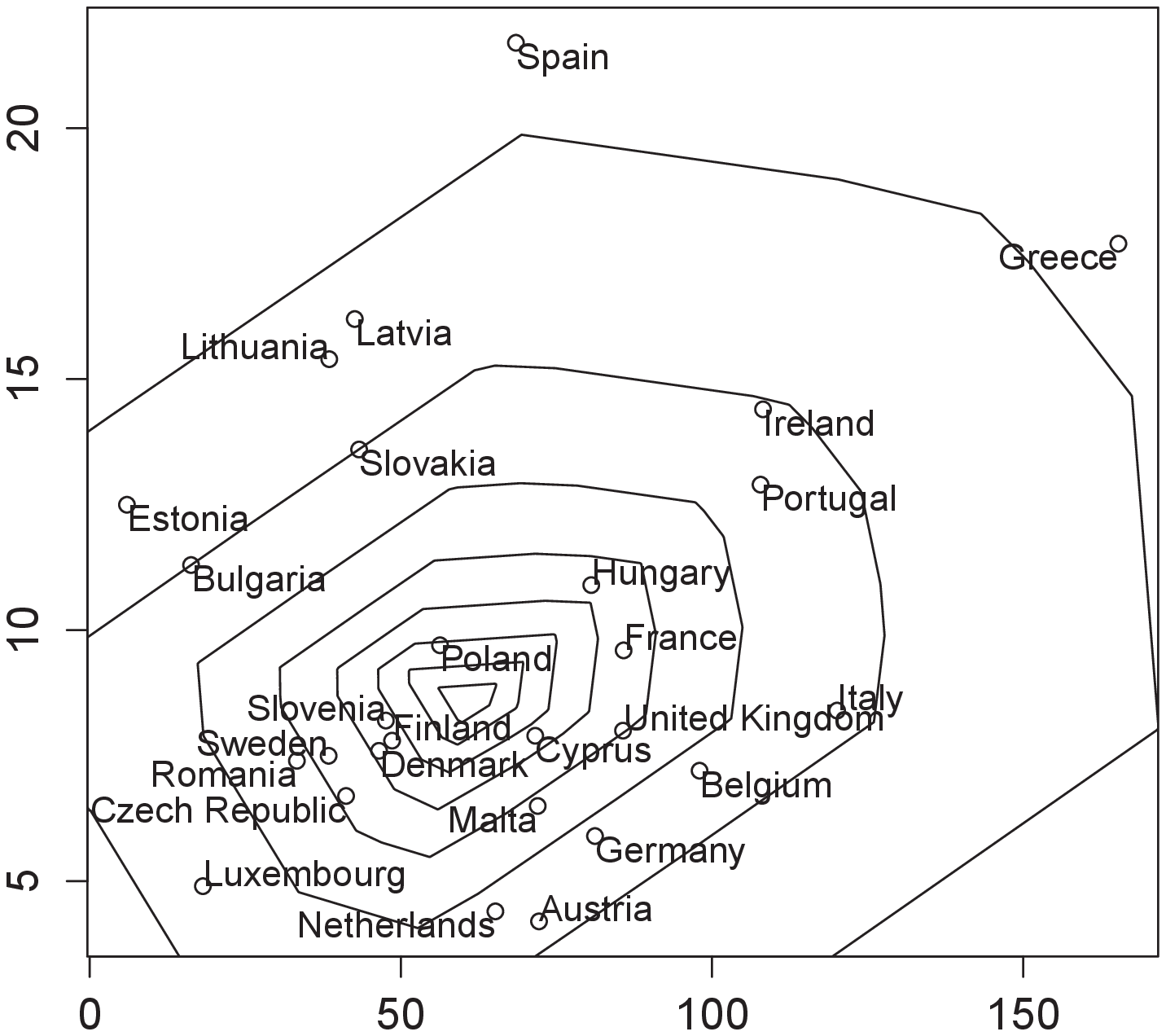}
\quad \includegraphics[scale=0.39]{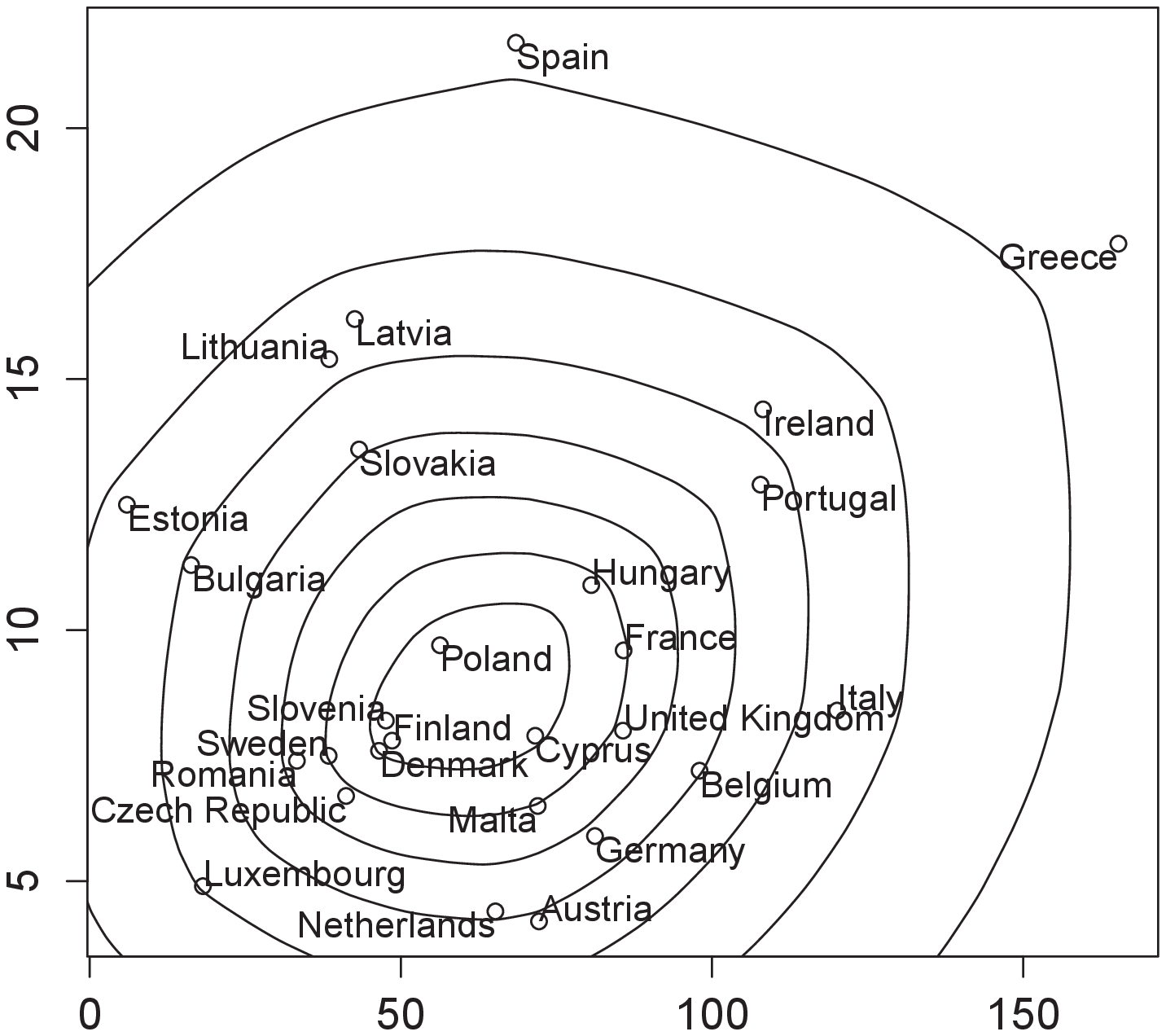}
\caption{Governmental debt and unemployment rate; projection depth regions (left), Oja regions (right); both with $\alpha = 0.1 (0.1),\dots, 0.9$.}
\label{fig2}
\end{figure}

{\bf Projection depth.}
The {\em projection depth} has been proposed in
\citet{ZuoS00a}:
\begin{equation}\label{projdepth}
\ddd{z}{X}{proj} =
         \left(1 + \sup_{p\in S^{d-1}}
         \frac {|\langle p,z \rangle
                - \med(\langle p,X \rangle)|}
               {\mbox{Dmed}(\langle p,X \rangle)}
         \right)^{-1} \,,
\end{equation}
where $S^{d-1}$ denotes the unit sphere in $\IR^d$,
$\langle p,z \rangle$ is the inner product (that is the projection of $z$ on the line $\{\lambda p : \lambda\in \IR\}$),
$\med(U)$ is the usual median of a univariate random variable $U$,
and $\mbox{Dmed}(U)=\med(|U-\med(U)|$ is the median absolute
deviation from the median.
The projection depth satisfies \textbf{D1} to \textbf{D5} and \textbf{D4con}. It has good properties, which are discussed in
detail by \citet{ZuoS00a}. For breakdown properties of the employed
location and scatter statistics, see \citet{Zuo00}.

{\bf Oja depth.} The {Oja depth} is not based on distances, but on average volumes of simplices
that have vertices from the data  \citep{ZuoS00a}:
%% see Figure \ref{fig2}.
$$\ddd{z}{X}{Oja} =\left( 1 + \frac
{\E \left(\V_d(\co\{z,X_{1},\ldots ,X_{d}\})\right)}
{\sqrt{\mbox{det\,}\Sigma_X}}
\right)^{-1}  ,
$$
where $X_1, \dots, X_d$ are random vectors independently distributed as $P$, $\co$ denotes the convex hull, $V_d$ the $d$-dimensional volume, and $S_{X}$ is defined as above. In particular, we can choose $D_{X}= {\boldmath \Sigma}_{X}$.
The Oja depth satisfies \textbf{D1} to \textbf{D5}.
It  is continuous on $z$ and maximum at the Oja median \citep{Oja83}, which is not
unique.
%%; see also Chapter \ref{Oja}.
The Oja depth determines the distribution uniquely among those measures
which have compact support of full dimension.

Figure \ref{fig2} contrasts the projection depth regions with the Oja regions for our debt-unemployment data.
The regions have different shapes, but agree in making Spain and Greece the most outlying countries.

\subsection{Weighted mean depths}\label{subsec3.2}

A large and flexible class of depth statistics corresponds to so called weighted-mean central regions, shortly WM regions \citep{DyckerhoffM11,DyckerhoffM10a}. These are convex compacts in $\IR^d$, whose support function is a weighted mean of order statistics, that is, an L-statistic.
Recall that a convex compact $K\subset\IR^d$ is uniquely determined by its support function $h_K$,
%% (see, e.g., \citet{Rockafellar70}),
\[
h_K(p)=\max\left\{p^\prime x\,:\,x\in K\,\right\}, \quad p\in S^{d-1}\,.
\]
To define the WM $\alpha$-region of an empirical distribution on $x^1, x^2, \dots,$ $ x^n$, we construct its support function as follows:
For $p\in S^{d-1}$, consider the line $\{\lambda p \in \IR^d : \lambda\in \IR\}$. By projecting the data on this line a linear ordering is obtained,
\begin{equation}\label{p-order}
	p^{\prime}x^{\pi_p(1)}\le p^{\prime}x^{\pi_p(2)}\le \dots \le p^{\prime}x^{\pi_p(n)}\,,
\end{equation}
and, by this, a permutation $\pi_p$ of the indices $1,2,\dots, n$.
%%Note that, if no equalities arise in (\ref{p-order}), the permutation $\pi_p$ is unique, otherwise a class $\Pi_p$ of several permutations is generated. The %%set of directions $p$ at which $\pi_p$ is not unique will be denoted $H(x^1,\dots,x^n)$,
%%\[
%%H(x^1,\dots,x^n)=\left\{p\in S^{d-1}\,:\,\text{there are $i\ne j$ such that
%%$p^\prime x^i=p^\prime x^j$}\,\right\}\,.
%%\]
Consider weights $w_{j,\alpha}$ for $j\in\{1,2,\dots,n\}$ and $\alpha\in [0,1]$ that satisfy
the following restrictions (i) to (iii):
\begin{enumerate}[(i)]
\item $\sum_{j=1}^n w_{j,\alpha}=1$, $w_{j,\alpha}\ge 0$ \ for all $j$ and $\alpha$\,.
\item $w_{j,\alpha}$ increases in $j$ \ for all $\alpha$\,.
\item $\alpha<\beta$ \; implies \quad
$\sum_{j=1}^kw_{j,\alpha}\le \sum_{j=1}^kw_{j,\beta}\,,\quad k=1,\dots,n\,.$
\end{enumerate}
Then, as it has been shown in \cite{DyckerhoffM11}, the function $h_{D_\alpha(x^1,\dots,x^n)}$,
\begin{equation}\label{defsuppWMT}
h_{D_\alpha(x^1,\dots,x^n)}(p)=\sum_{j=1}^nw_{j,\alpha}p^\prime x^{\pi_p(j)}\,, \ p\in S^{d-1}\,,
\end{equation}
is the support function of a convex body $D_\alpha=D_\alpha(x^1,\dots,x^n)$, and $D_\alpha \subset D_\beta$ holds whenever $\alpha> \beta$.
Now we are ready to see the general definition of a family of WM regions.
\begin{definition}
Given a weight vector $w_\alpha=w_{1,\alpha},\dots w_{n,\alpha}$ that satisfies the restrictions
$(i)$ to $(iii)$, the convex compact $D_\alpha=D_\alpha(x^1,\dots,x^n)$ having support function (\ref{defsuppWMT})
is named the \emph{WM region} of $x^1,\dots,x^n$ at level $\alpha$\,, $\alpha \in [0,1]$.
The corresponding depth (\ref{eqregiondepth}) is the \emph{WM depth} with
weights $w_\alpha$, $\alpha\in [0,1]$.
\end{definition}
It follows that the WM depth satisfies the restrictions \textbf{D1} to \textbf{D5} and \textbf{D4con}.
Moreover, it holds
\begin{equation}
\label{eqregconv}
D_\alpha(x^1,\dots,x^n)=conv\left\{\sum_{j=1}^nw_{j,\alpha}x^{\pi(j)}\, : \,
\text{$\pi$ permutation of $\{1,\dots,n\}$}\,\right\}\,.
\end{equation}
%%and the set of extreme points of $D_\alpha$ is given by
%%\begin{equation}
%%\label{eqext}
%%Ext\bigl(D_\alpha(x^1,\dots,x^n)\bigr)=\left\{\sum_{j=1}^nw_{j,\alpha}x^{\pi_p(j)}\,:\,
%%p\in S^{d-1}\setminus H(x^1,\dots,x^n)\,\right\}\,.
%%\end{equation}

This explains the name by stating that a WM region is the convex hull of weighted means of the data.
Consequently, outside the convex hull of the data the WM depth vanishes.
WM depths are useful statistical tools as their central regions have attractive analytical and computational properties. Sample WM regions are consistent estimators for the WM region of the underlying probability. Besides being \textit{continuous} in the distribution and in $\alpha$, WM regions  are \textit{subadditive}, that is,
\[ D_\alpha(x^1+y^1,\dots,x^n+y^n) \subset  D_\alpha(x^1,\dots,x^n) \oplus D_\alpha(y^1,\dots,y^n)\,,
\]
and \textit{monotone}:
If $x^i\le y^i$  holds for all $i$ (in the componentwise ordering of $\IR^d$), then
\begin{eqnarray*}
&& D_\alpha(y^1,\dots,y^n) \subset D_\alpha(x^1,\dots,x^n) \oplus \IR^d_+ \quad \text{and} \\
&& D_\alpha(x^1,\dots,x^n) \subset D_\alpha(y^1,\dots,y^n) \oplus \IR^d_- \,,
\end{eqnarray*}
where $\oplus$ signifies the Minkowski sum of sets.

Depending on the choice of the weights $w_{j,\alpha}$ different notions of data depths are obtained.
%%They include the zonoid depth \citep{KoshevoyM97b}, the expected convex hull (ECH$^*$) depth \citep{Cascos07}, the geometrical depth., and others.
For a detailed discussion of these and other special WM depths and central regions, the reader is referred to \cite{DyckerhoffM11,DyckerhoffM10a}.

\begin{figure}
\includegraphics[scale=0.39]{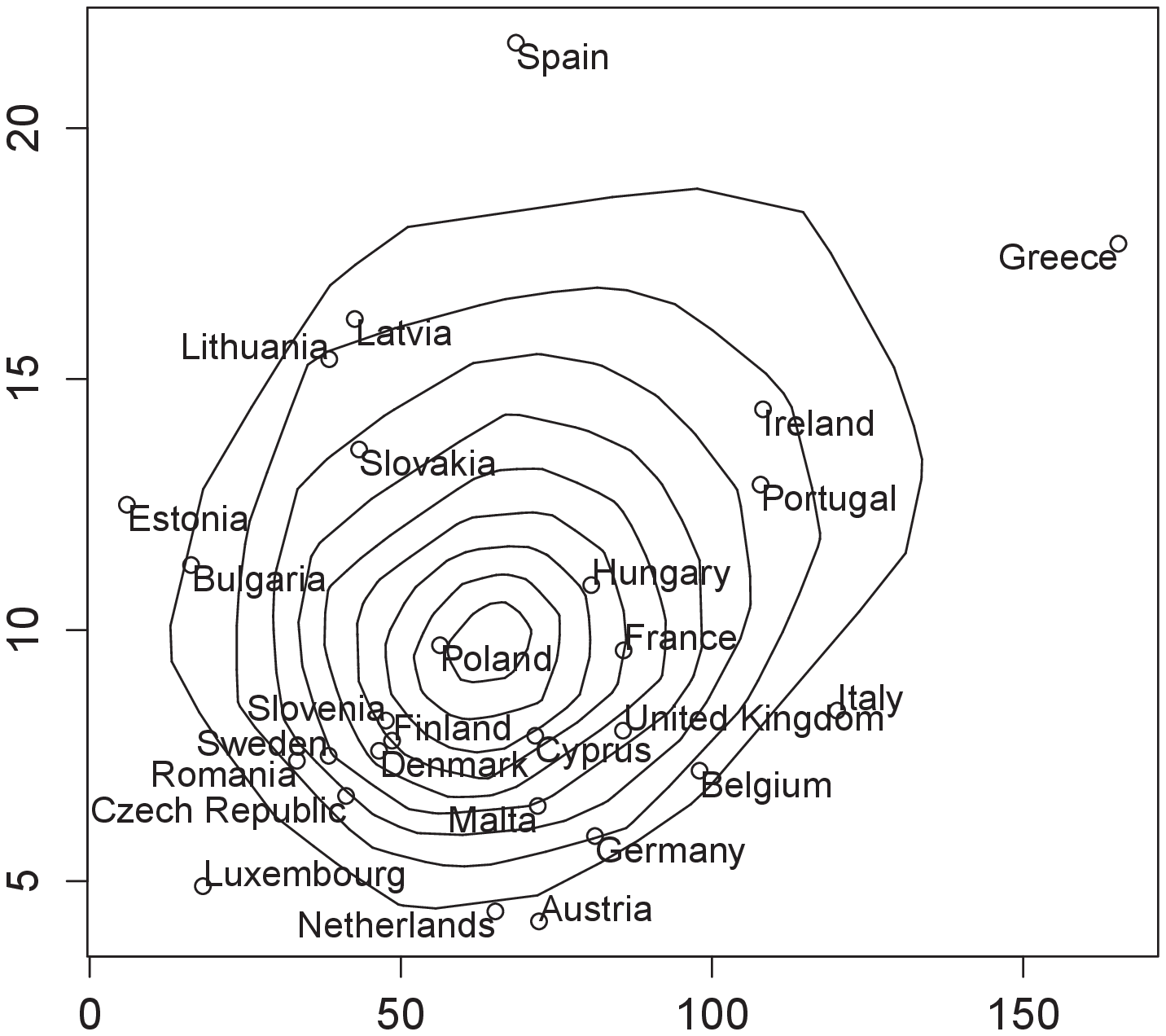}
\quad \includegraphics[scale=0.39]{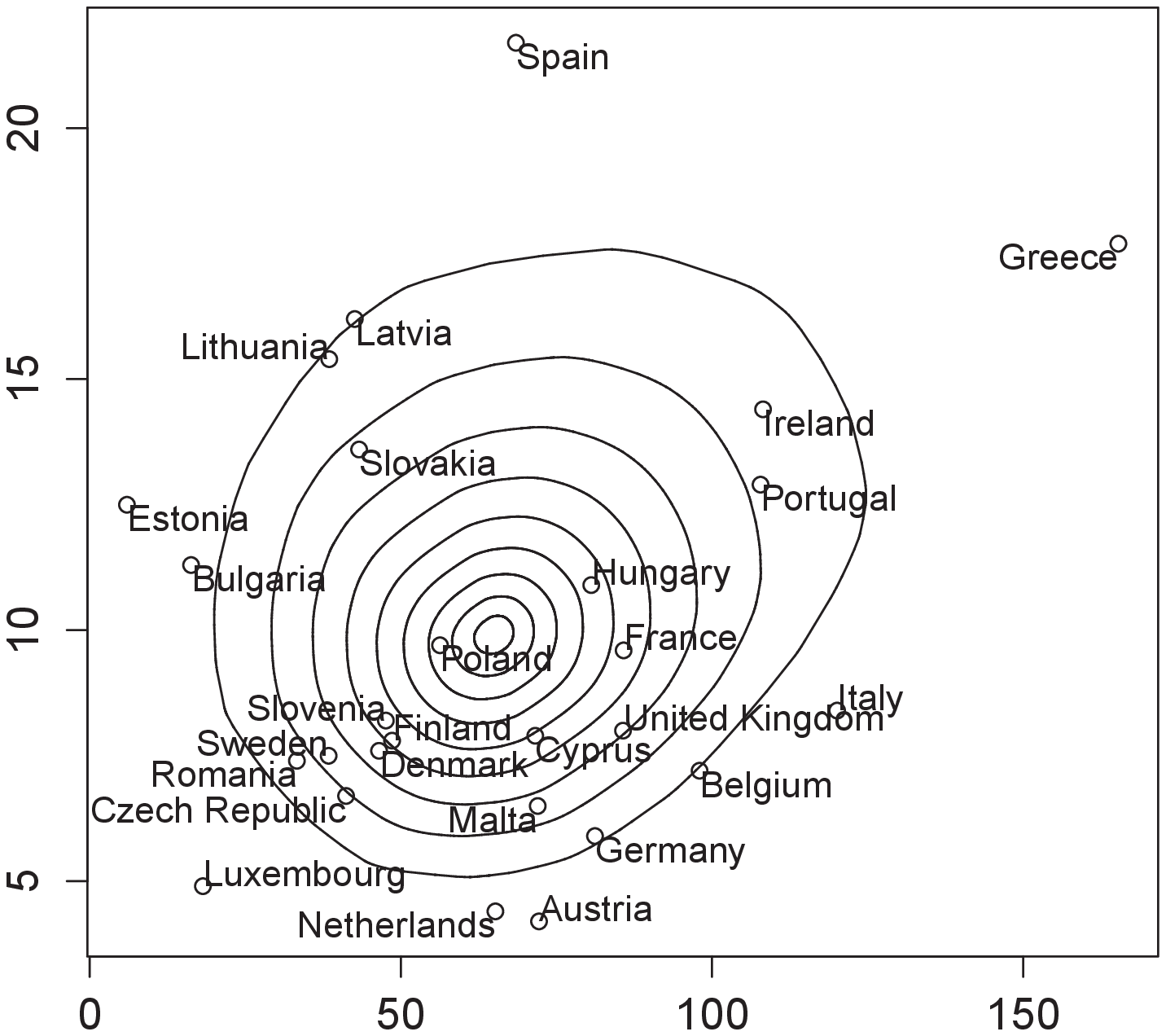}
\caption{Governmental debt and unemployment rate; zonoid regions (left), ECH$^*$ regions (right); both with $\alpha = 0.1 (0.1),\dots, 0.9$.}
\label{fig3}
\end{figure}

\textbf{Zonoid depth.}
%%Zonoid regions have been introduced by \citet{KoshevoyM97b}.
For an empirical distribution $P$ on $x^1,\dots, x^n$ and $0<\alpha\le1$ define the zonoid region \citep{KoshevoyM97b}
\[
D^{zon}_\alpha(P)=\left\{\sum_{i=1}^n\lambda_ix^i\,:\,
0\le \lambda_i\le\frac{1}{n\alpha}\,,\, \sum_{i=1}^n\lambda_i=1\right\}\,.
\]
See Figure \ref{fig3}.
The corresponding support function (\ref{defsuppWMT})
%%\[
%%h_{D^{zon}_\alpha}(p)=\sum_{j=1}^nw_{j,\alpha}p'x^{\pi_p(j)}\,,
%%\]
employs the weights
\begin{equation}
\label{eqzonweights}
w_{j,\alpha}=\left\{\begin{array}{cl}
0\,&\text{if $j<n-\lfloor n\alpha\rfloor$}\,,\\[1ex]
\frac{n\alpha-\lfloor n\alpha\rfloor}{n\alpha}\,&\text{if $j=n-\lfloor n\alpha\rfloor$}\,,\\[1ex]
\frac{1}{n\alpha}\,&\text{if $j>n-\lfloor n\alpha\rfloor$}\,.
\end{array}\right.
\end{equation}
Many properties of zonoid regions and the zonoid depth $D^{zon}(z|X)$ are discussed in \cite{Mosler02a}.
The zonoid depth lift equals the so called lift zonoid, which fully characterizes the distribution. Therefore the zonoid depth generates an antisymmetric depth order (\ref{eqdepthordering}) and a probability metric (\ref{eqdepthnorm}). Zonoid regions are not only invariant to affine, but to general linear transformations; specifically any marginal projection of a zonoid region is the zonoid region of the marginal distribution.
The zonoid depth is continuous on $z$ as well as $P$.

\textbf{Expected convex hull depth.}
Another important notion of WMT depth is that of \emph{expected convex hull (ECH*)} depth \citep{Cascos07}.
%%; see Figure \ref{fig3}.
Its central region $D_\alpha$ (see Figure \ref{fig3}) has a support function
%%\[
%%h_{ECH^*_\alpha}(p)=\sum_{j=1}^n w_{j,\alpha}\,p^\prime x^{\pi_p(j)}\,,
%%\]
with weights
\begin{equation}\label{weightsECHstar}
w_{j,\alpha}=\frac{j^{1/\alpha}-(j-1)^{1/\alpha}}{n^{1/\alpha}}\,.
\end{equation}
Figure 3 depicts zonoid and ECH$^*$ regions for our data. We see that the zonoid regions are somewhat angular while the ECH$^*$ regions appear to be smoother;
this corresponds, when calculating such regions in higher dimensions, to a considerably higher computation load of ECH$^*$.

\textbf{Geometrical depth.} The weights
\[
w_{j,\alpha}=\left\{
\begin{array}{cl}
\frac{1-\alpha}{1-\alpha^n}\,\alpha^{n-j} &\text{if $0<\alpha<1$\,,}\\
 0                                        &\text{if $\alpha=1$\,,}
\end{array}\right.
\]
yield another class of WM regions. The respective depth is the \emph{geometrically weighted mean depth} \citep{DyckerhoffM11}.

\subsection{Depths based on halfspaces and simplices}\label{subsec3.3}

The third approach concerns no distances or volumes, but the combinatorics of halfspaces and simplices only. In this it is independent of the metric structure of $\IR^d$. While depths that are based on distances or weighted means may be addressed as \emph{metric depths}, the following ones will be mentioned as \emph{combinatorial depths}. They remain constant, as long as the compartment structure of the data does not change.
By this, they are very robust against \emph{location outliers}. Outside the convex support $\conv(X)$ of the distribution every combinatorial depth attains its minimal value, which is zero.

\begin{figure}
\includegraphics[scale=0.39]{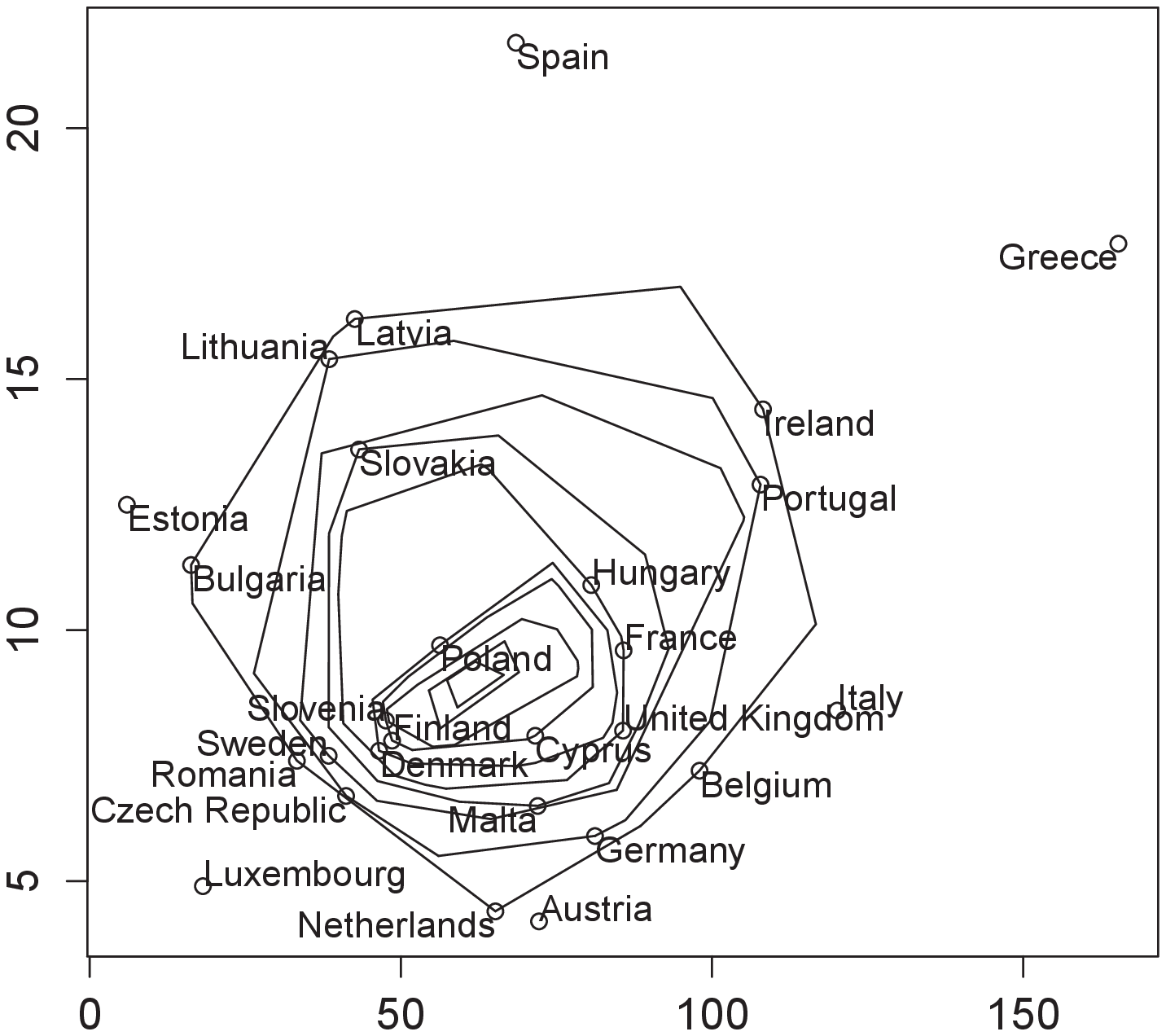}
\quad \includegraphics[scale=0.39]{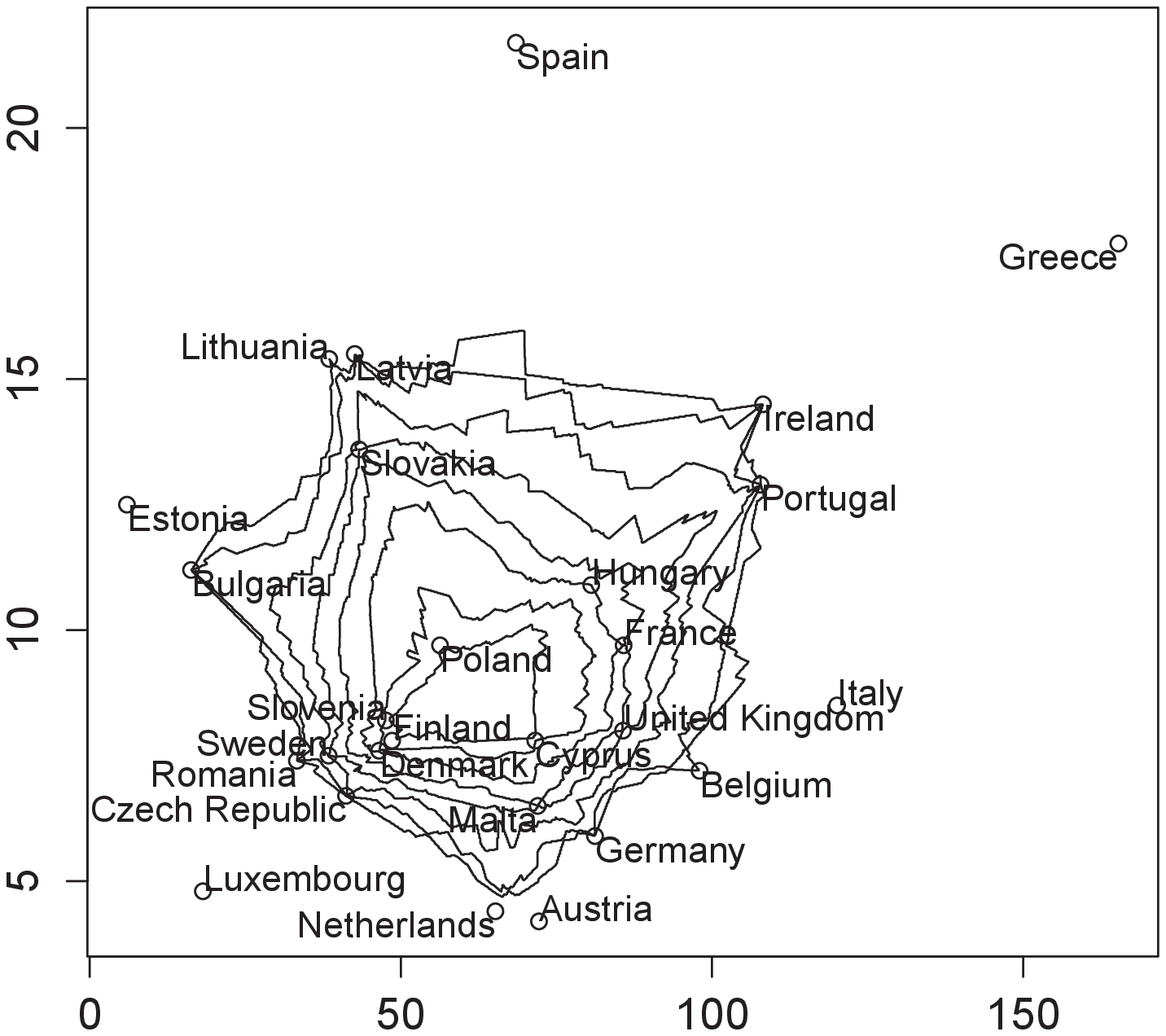}
\caption{Governmental debt and unemployment rate; Tukey regions (left) with $\alpha=\frac 2{27} (\frac 1{27}),\dots, \frac 11{27}$, simplicial regions (right) with $\alpha= 0.25, 0.3 (0.1),\dots, 0.9$.}
\label{fig4}
\end{figure}

\textbf{Location depth.}
Consider the population version of the \emph{location depth},
\begin{equation}\label{Hdepth}
\ddd{z}{X}{loc} =
\inf \{P(H) :  H \;\; \mbox{is a closed halfspace,} \;\; z \in H\}\,.
\end{equation}
The depth is also known as \emph{halfspace} or \emph{Tukey depth}, its central regions as \emph{Tukey regions}.
%%The sample version is in the same way defined on the empirical distribution
%%of the data.
The location depth is affine invariant (\textbf{D1}, \textbf{D2}).
Its central regions are convex (\textbf{D4con}) and closed (\textbf{D5}); see Figure \ref{fig4}.
The maximum value of the location depth is smaller or equal to $1$ depending on the distribution.
The set of all such points is mentioned as the {\em halfspace median set}
and each of its elements as a {\em Tukey median}\/ \citep{Tukey75}.

If $X$ has an \emph{angular symmetric} distribution, the location depth attains its  maximum at the center
and the center is a Tukey median; this strengthens Proposition \ref{prop:symmetric}.
(A distribution is called \emph{angular} (= \emph{halfspace}) \emph{symmetric} about $z^*$  if $P(X \in H)\ge 1/2$ for every closed halfspace H having $z^*$ on the boundary; equivalently, if $(X-z^*)/||X-z^*||$ is centrally symmetric with the convention $0/0=0$.)

If $X$ has a Lebesgue-density, the location depth depends continuously
on $z$; otherwise
the dependence on $z$ is noncontinuous and there can be more
than one point where the maximum is attained.
As a function of $P$ the location depth is obviously
noncontinuous. It determines the distribution in a unique way if
the distribution is either  discrete
\citep{StruyfR99,Koshevoy99c} or continuous with compact support.
The location depth of a sample from $P$
converges almost surely to the location depth of $P$
\citep{DonohoG92}.
%%The same holds for the depth contours
%%\citep{MasseT94,HeW97,ZuoS00d} if the distribution is elliptic.
The next depth notion involves simplices in $\IR^d$.

\textbf{Simplicial depth.}
\citet{Liu90} defines the {\em simplicial depth}
as follows:
\begin{equation}\label{Simdepth}
\ddd{z}{X}{sim} = P\left(z \in
%{\it int} \,
\co (\{ X_{1},\ldots ,X_{{d+1}} \})\right) ,
\end{equation}
where $X_{1},\ldots ,X_{{d+1}}$ are i.i.d.\ by $P$.
The sample version reads as
\begin{equation}\label{sSimdepth}
\dddd{z}{x}{sim}=
\frac 1{{n\choose{d+1}}}\
\#\Bigl\{ \{i_1,\ldots ,i_{d+1}\}\,:\,z\in
%{\it int} \,
\co (\{ x^{i_1},\ldots ,x^{i_{d+1}} \})\Bigr\}.
\end{equation}
The simplicial depth is affine invariant (\textbf{D1}, \textbf{D2}).
%%Again, it vanishes outside the $co(P)$ and
Its maximum is less or equal to $1$, depending on the distribution.
In general, the point of maximum simplicial depth is not unique;
the {\em simplicial median}\/ is defined as the gravity center of
these points.
The sample simplicial depth converges almost surely uniformly in
$z$ to its population version \citep{Liu90,Duembgen92}.
The simplicial depth has positive breakdown \citep{Chen95}.

If the distribution is Lebesgue-continuous, the simplicial depth behaves
well: It varies continuously on $z$ \citep[Th.~2]{Liu90},
is maximum at a center of angular symmetry,
and decreases monotonously from a deepest point (\textbf{D4}).
The {\em simplicial central regions} of a Lebesgue-continuous
distribution are connected and compact \citep{Liu90}.

However, if the distribution is discrete, each of these properties can
fail; for counterexamples see, e.g., \citet{ZuoS00a}.
The {simplicial depth} characterizes an empirical measure
if the supporting points are in \emph{general position}, that is,
if no more than $d$ of the points lie on the same hyperplane.

As Figure \ref{fig4} demonstrates, Tukey regions are convex while simplicial regions are only starshaped. The Figure illustrates also that these notions are rather insensitive to outlying data: both do not reflect {\em how far}\/ Greece and Spain are from the center. Whether, in an application, this kind of robustness is an advantage or not, depends on the problem and data at hand.

Other well known combinatorial data depths are the {\em majority depth} \citep{LiuS93}
and the \emph{convex-hull peeling depth} \citep{Barnett76,DonohoG92}. However the latter possesses no population version.

\section{Functional data depth}\label{sec4}
The analysis of functional data has become a practically important branch of statistics; see \cite{RamsayS05}.
Consider a space $E$ of functions $[0,1]\to \IR$ with the supremum norm.
Like a multivariate data depth, a functional data depth is a real-valued functional that indicates how `deep' a function $z\in E$ is located in a given finite cloud of functions $\in E$. Let $E'$ denote the set of continuous linear functionals $E\to \mathbb{R}$, and ${E'}^d$ the $d$-fold Cartesian product of $E'$.
Here, following \cite{MoslerP12}, functional depths of a general form (\ref{defFD}) are presented.
%%This definition includes several notions of functional data depth that are known from the literature as well as many new ones.
Some alternative approaches will be addressed below.

\textbf{$\Phi$-depth.}
For $z\in E$ and an empirical distribution $X$ on $x^1, \dots, x^n\in E$, define a \emph{functional data depth} by
\begin{equation}\label{defFD}
 D(z|X)= \inf_{\varphi\in \Phi} D^d(\varphi(z)|\varphi(X))\,,
 \end{equation}
where $D^d$ is a $d$-variate data depth satisfying {\rm \textbf{D1}} to {\rm \textbf{D5}},  $\Phi\subset {E'}^d$, and $\varphi(X)$ is the empirical distribution on $\varphi(x^1), \dots, \varphi(x^n)$.  $D$ is called a \emph{$\Phi$-depth}. A population version is similarly defined.

Each $\varphi$ in this definition may be regarded as a particular `aspect' we are interested in and which is represented in $d$-dimensional space. The
depth of $z$ is given as the smallest multivariate depth of $z$ under all these aspects. It implies that all aspects are equally relevant so that the depth of $z$ cannot be larger than its depth under any aspect.

As the $d$-variate depth $D^d$ has maximum not greater than $1$, the functional data depth $D$ is bounded above by 1. At every point $z^*$ of \textit{maximal $D$-depth} it holds $D(z^*|X)\le 1$. The bound is attained with equality, $D(z^*|X)= 1$, iff
$D^d(\varphi(z^*)|\varphi(X))=1$ holds for all $\varphi\in \Phi$, that is, iff
\begin{equation}\label{deepest}
z^* \in \bigcap_{\varphi \in \Phi} \varphi^{-1}(D_1^d(\varphi(X)))\,.
\end{equation}

A $\Phi$-depth (\ref{defFD}) always satisfies {\rm \textbf {D1}}, {\rm \textbf {D2sca}}, {\rm \textbf {D4}}, and {\rm \textbf {D5}}.

   It satisfies {\rm \textbf {D3}} if for every sequence $(z^i)$ with $||z^i||\to \infty$  exists a  $\varphi$ in $\Phi$ such that $\varphi(z^i)\to \infty$\,. (For some special notions of functional data depth this postulate has to be properly adapted.)

   {\rm \textbf{D4con}} is met if {\rm \textbf{D4con}} holds for the underlying $d$-variate depth.

 We now proceed with specifying the set $\Phi$ of functionals and the multivariate depth $D^k$ in (\ref{defFD}). While many features of the functional data depth (\ref{defFD}) resemble those of a multivariate depth, an important difference must be pointed out:
In a general Banach space the unit ball $B$ is not compact, and properties \textbf{D3} and  \textbf{D5} do not imply that the level sets of a functional data depth are compact.
So, to obtain a meaningful notion of functional data depth of type (\ref{defFD}) one has to carefully choose a set of functions $\Phi$ which is not too large. On the other hand, $\Phi$ should not be too small, in order to extract sufficient information from the data.

\textbf{Graph depths.}
%%Consider $E=C\left(J;\mathbb{R}^d\right)$ with norm $||\cdot||_\infty$ and let
For $x\in E$ denote $x(t)=(x_1(t)\ldots ,x_d(t))$ and consider
\begin{equation}\label{bandPhi}
\Phi= \{\varphi^t:E\to \IR^d \,:\, \varphi^t (x)=(x_1(t)\ldots ,x_d(t)),\, t\in T\}
\end{equation}
for some $T\subset [0,1]$, which may be a subinterval or a finite set. For $D^d$ use any multivariate depth that satisfies \textbf{D1} to \textbf{D5}. This results in the \emph{graph depth}
\begin{equation}\label{defBand}
 GD(z|x^1,\dots, x^n)= \inf_{t\in T} D^d(z(t)|x^1(t),\dots, x^n(t))\,.
\end{equation}

In particular, with the univariate halfspace depth, $d=1$ and $T=J$ we obtain the \emph{halfgraph depth} \citep{LopezPR05}.
Also, with the univariate simplicial depth the \emph{band depth} \citep{LopezPR09} is obtained, but this, in general, violates monotonicity \textbf{D4}.

\textbf{Grid depths.}
We choose a finite number of points in $J$, $t_1,\dots, t_k$, and evaluate a function $z\in E$ at these points. Notate $\underline t=(t_1,\dots, t_k)$ and $z(\underline t)=(z_1(\underline t), \dots, z_d(\underline t))^\textsf{T}$.  That is, in place of the function $z$ the $k\times d$ matrix
$z^{(k)}$ is considered. A \emph{grid depth} $RD$ is defined by (\ref{defFD}) with the following $\Phi$,
\begin{equation}
\Phi= \{\varphi^r : \varphi^r(z)= (\langle r,z_1(\underline t)\rangle, \dots, \langle r,z_d(\underline t)\rangle), r\in S^{k-1}\}\,,
\end{equation}
which yields
\begin{equation}\label{defGrid}
 RD(z|x^1,\dots, x^n)= \inf_{r\in S^{k-1}} D^d(\langle r, z(\underline t)\rangle| \langle r, x^1(\underline t)\rangle,\dots ,\langle r, x^n(\underline t)\rangle)   \,.
\end{equation}

A slight extension of the $\Phi$-depth is the \emph{principal components depth} \citep{MoslerP12}.
However, certain approaches from the literature are no $\Phi$-depths.
These are mainly of two types. The first type employs \emph{random projections} of the data:
\cite{CuestaAN08b} define the depth of a function as the univariate depth of the function values taken at a randomly chosen argument $t$.
\cite{CuevasFF07} also employ a random projection method.
The other type uses average univariate depths.
%%Compared to this our definition may be mentioned as a `uniform' depth.
\cite{FraimanM01} calculate the univariate depths of the values of a function and integrate them over the whole interval; this results in kind of `average' depth.
\cite{ClaeskensHS12} introduce a multivariate ($d\ge 1$) functional data depth, where they similarly compute a weighted average depth. The weight at a point reflects the variability of the function values at this point (more precisely: is proportional to the volume of a central region at the point).
%%These notions satisfy the above basic postulates or proper modifications of them; but a detailed analysis of them is beyond the scope of this paper.

\section{Computation of depths and central regions}\label{sec5}
The moment Mahalanobis depth and its elliptical central regions are obtained in any dimension by calculating the mean and the sample covariance matrix, while robust Mahalanobis depths and regions are determined with the R-procedures ``cov.mcd'' and ``cov.mve''.
In dimension $d=2$, the central regions of many depth notions can be exactly calculated by following a \emph{circular sequence} \citep{Edelsbrunner87}.
The R-package ``depth'' computes the exact location ($d=2,3$) and simplicial ($d=2$) depths, as well as the Oja depth and an approximative location depth for any dimension.
An exact algorithm for the location depth in any dimension is developed in \cite{LiuZ12}.
\cite{CuestaAN08a} propose to calculate instead the \emph{random Tukey depth}, which is the minimum univariate location depth of univariate projections in a number of randomly chosen directions.
With the algorithm of \cite{PaindaveineS12}, Tukey regions are obtained, $d\ge 2$.
The bivariate projection depth is computed by the R-package "ExPD2D"; for the respective regions, see \cite{LiuZW11}.
The zonoid depth can be  efficiently determined in any dimension \citep{DyckerhoffKM96}.
An R-package (``WMTregions'') exists for the exact calculation of zonoid and general WM regions; see
\cite{MoslerLB09, BazovkinM12}. The R-package ``rainbow'' calculates several functional data depths.

\section{Conclusions}\label{sec6}
Depth statistics have been used in numerous and diverse tasks of which we can mention a few only. \cite{LiuPS99} provide an introduction to some of them.
In descriptive multivariate analysis, depth functions and central regions visualize the data regarding location, scale and shape. By bagplots and sunburst plots outliers can be identified and treated in an interactive way. In $k$-class supervised classification, each - possibly high-dimensional - data point is represented in $[0,1]^k$ by its values of depth in the $k$ given classes, and classification is done in $[0,1]^k$,
Functions of depth statistics include depth-weighted statistical functionals, such as
${\int_{\IR^d} x\, w(D(x|P)) dP}/{\int_{\IR^d}  w(D(x|P)) dP}$ for location.
In inference, tests for goodness of fit and homogeneity regarding location, scale and symmetry are based on depth statistics; see, e.g.\ \cite{Dyckerhoff02, LeyP11}. Applications include such diverse fields as statistical control \citep{LiuS93}, measurement of risk \citep{CascosM07}, and robust linear programming \citep{BazovkinM11}.
Functional data depth is applied to similar tasks in description, classification and testing; see e.g.\ \cite{LopezPR09, CuevasFF07}.

This survey has covered the fundamentals of depth statistics for $d$-variate and functional data.
Several special depth functions in $\IR^d$ have been presented, metric and combinatorial ones, with a focus on the recent class of WM depths.
For functional data, depths of infimum type have been discussed.
%%Other functional data depths employ averages, that is integrals of uni- or multivariate depths over the domain; see e.g.\ \cite{CuevasF09}, %%\cite{ClaeskensHS12}.
Of course, such a survey is necessarily incomplete and biased by the preferences of the author. Of the many applications of depth in the literature only a few have been touched, and important theoretical extensions like regression depth \citep{RousseeuwH96}, depth calculus \citep{Mizera02},
location-scale depth \citep{MizeraM04}, and likelihood depth \citep{Mueller05}
have been completely omitted.

Most important for the selection of a depth statistic in applications are the questions of computability and - depending on the data situation - robustness. Mahalanobis depth is solely based on estimates of the mean vector and the covariance matrix. In its classical form with moment estimates Mahalanobis depth is efficiently calculated but highly non-robust, while with estimates like the minimum volume ellipsoid it becomes more robust. However, since it is constant on ellipsoids around the center, Mahalanobis depth cannot reflect possible asymmetries of the data. Zonoid depth can be efficiently calculated, also in larger dimensions, but has the drawback that the deepest point is always the mean, which makes the depth non-robust. So, if robustness is an issue, the zonoid depth has to be combined with a proper preprocessing of the data to identify possible outliers. The location depth is, by construction, very robust but expensive when exactly computed in dimensions more than two.
As an efficient approach the random Tukey depth yields an upper bound on the location depth, where the number of directions has to be somehow chosen.

A depth statistics measures the centrality of a point in the data. Besides ordering the data it provides numerical values that, with some depth notions, have an obvious meaning; so with the location depth and all WM depths. With other depths, in particular those based on distances, the outlyingness function has a direct interpretation.

%%More research is needed in three fields: Finding multivariate (non-convex, local) depths for the analysis of multimodal distributions,
%%constructing and applying depth statistics in more general spaces, developing efficient exact and approximate algorithms to calculate depths and central %%regions in higher dimensions.

%%Other qualified candidates to be calculated in higher dimensions are - albeit being only mirror symmetric - $L^p$-depths.

\bibliographystyle{annals}

\bibliography{DEPTHSTAT}

\end{document}